\documentclass[12pt]{article}
\usepackage{amssymb}
\newtheorem{Thm}{Theorem}

\newtheorem{thm}{Theorem}[section]
\newtheorem{dfn}{Definition}[section]
\newtheorem{prop}{Proposition}[section]
\newtheorem{lem}{Lemma}[section]
\newtheorem{rem}{Remark}[section]

\def\A{{\cal{A}}}
\def\B{{\cal{B}}}

\def\C{{\mathbb{C}}}

\def\Z{{\mathbb{Z}}}
\def\ka{{\kappa}}

\def\la{{\langle}}
\def\ra{{\rangle}}
\def\lla{{\langle \! \langle}} 
\def\rra{{\rangle \! \rangle}}
\title{Nichols-Woronowicz model of coinvariant algebra of 
complex reflection groups} 
\date{}
\author{Anatol N. Kirillov and Toshiaki Maeno}
\begin{document}
\maketitle
\begin{abstract}
We give a model of the coinvariant algebra of the complex reflection 
groups as a subalgebra of a braided Hopf algebra called Nichols-Woronowicz 
algebra. 
\end{abstract}
\section*{Introduction}
Let $V$ be a finite dimensional complex vector space. A finite subgroup 
$G \subset GL(V)$ is called {\it a complex reflection group}, if $G$ can be 
generated by the set of pseudo-reflections, i.e., transformations that fix a complex
hyperplane in $V$ pointwise. Any real reflection group becomes a complex 
reflection group if one extends the scalars from $\mathbb R$ to $\mathbb C$. 
In particular all Coxeter groups give examples of complex reflection groups. 
We refer the reader to \cite{BMR} for general background of the theory of 
complex reflection groups. Below we recall a few facts about real and complex 
reflection groups which appeared to be a motivation for our paper. 

In 1954, G. C. Shephard and J. A. Todd \cite{ST} had obtained a complete 
classification of finite irreducible complex reflection groups. They found 
that there exist an infinite family of irreducible complex reflection groups
$G(e,p,n)$ depending on three positive integer parameters (with $p$ dividing 
$e$), and 34 exceptional groups $G_4, \ldots,G_{37}.$ The group 
$G(e,p,n)$ has the order $e^nn!/p.$ It also has a normal abelian subgroup of 
order $e^n/p,$ and the corresponding quotient is the symmetric group 
on $n$ points. The family of groups $G(e,p,n)$ includes the cyclic group $C_{e/p}$ of 
order $e/p,$ namely, $C_{e/p}=G(e,p,1);$ the symmetric group on $n$ 
points $S_n=G(1,1,n);$ the Weyl groups of types $B_n,$ $C_n,$ and 
$D_n,$ namely, $W_{B_n}=W_{C_n}= G(2,1,n)$ and $W_{D_n}=G(2,2,n);$ and the dihedral groups 
$I_2(e)=G(e,e,2).$ 

The fundamental fact characterizing the finite complex reflection subgroups in
$GL(V)$ is the following theorem by G. C. Shephard and J. A. Todd.
\begin{Thm} {\rm (Shephard-Todd \cite{ST})} A subgroup $G \subset GL(V)$ is a finite 
complex reflection group 
if and only if the subring $P^{G}$ of the $G$-invariant elements in the symmetric 
algebra $P=S(V)$ of the space $V$ is generated by $n$ algebraically 
independent homogeneous elements. 
\end{Thm}
On the other side, in the case of Coxeter groups that form a part of the
complex reflection groups, there is a remarkable result by C. F. Dunkl which 
states that the algebra generated by the truncated Dunkl operators is 
isomorphic to the coinvariant algebra of the corresponding Coxeter group. 
An analogue of Dunkl operators for finite complex reflection groups have been 
introduced by C. F. Dunkl and E. M. Opdam \cite{DO}. So it seems an interesting 
problem to extend the result by C. F. Dunkl mentioned above, to the case of 
complex reflection groups.

It is well-known that the cohomology ring of a flag variety has a 
presentation as the coinvariant algebra of the corresponding Weyl group. 
Some combinatorial problems on the intersection theory over flag varieties 
can be formulated for the coinvariant algebra of a finite Coxeter group 
\cite{Hi1}. In view of Shephard and Todd's theorem, the coinvariant algebra 
of a finite complex reflection group gives a natural generalization of 
the framework where one can study problems related to the Schubert calculus, 
see e.g. \cite{To}. 

S. Fomin and the first author \cite{FK} have given a model of the cohomology 
ring of the flag variety of type $A$ as a commutative subalgebra in a certain 
noncommutative quadratic algebra. Their construction has applications to 
Pieri's formula, quantization and so on \cite{FK}, \cite{Po}. Similar 
construction for other root systems has been given in \cite{KM1}. Yu. Bazlov 
\cite{Ba} has realized the coinvariant algebra of a finite Coxeter group as a 
commutative subalgebra in a braided Hopf algebra, called the Nichols-Woronowicz algebra, 
to give a new mode of thought on the construction in \cite{FK}. The 
quantization operator on the Nichols-Woronowicz algebra and 
the model of the quantum cohomology ring of the flag 
varieties are given in \cite{KM2}. 

The Nichols-Woronowicz algebra $\B(M),$ which is called the Nichols algebra in 
\cite{AS}, associated to a braided vector space $M$ is 
a braided graded Hopf algebra characterized by the following condition which 
appeared originally in the work of W. D. Nichols \cite{Ni}: \\ 
(1) $\B^0(M)=\C,$ \\ 
(2) $\B^1(M)=M=\{ \textrm{primitive elements in $\B(M)$} \},$ \\ 
(3) $\B^1(M)$ generates $\B(M)$ as an algebra. \\ 
It is known that the algebra $\B(M)$ has an alternative definition 
as the braided analogue of the symmetric (or exterior) algebras introduced 
by S. L. Woronowicz \cite{Wo} for the study of differential calculus 
on quantum groups, see \cite{Sc}. 

In the present paper, we give a generalization of Bazlov's construction 
to the case of finite complex reflection groups. Having this aim in mind, we 
define the Yetter-Drinfeld module $M_{G}$ corresponding to a finite complex 
reflection group $G.$  It is similar to the case of finite Coxeter groups that 
a linear basis of the Yetter-Drinfeld module $M_{G}$ is parametrized by the set of 
pseudo-reflections in the group $G.$ However, the $G$-grading and the $G$-module 
structure on $M_{G}$ essentially depend on the properties of the 
hyperplane arrangement $\A=\A_{G}$ consisting of the reflection hyperplanes 
corresponding to the group $G,$ see 
Section 2 for details. With the Yetter-Drinfeld module $M_{G}$ in hand, 
the construction of the Woronowicz symmetrizers $\sigma_n$ and the corresponding 
Nichols-Woronowicz algebra ${\cal B}(M_G)$ is done in the standard manner, see \cite{Ba}, 
\cite{MS}. In Proposition 3.2 we compute the set of quadratic relations in the 
algebra ${\cal B}(M_G)$ in the case $G=G(e,1,n).$ In Section 4 we 
construct a realization of the coinvariant algebra $P_G$ of a finite complex 
reflection group $G$ as a commutative subalgebra in the corresponding 
Nichols-Woronowicz algebra ${\cal B}(M_G).$ The basis for our construction 
is the definition of the $\C$-linear map $\mu,$ see Definition 4.1. The map 
$\mu$ can be treated as "a truncated version" of the Dunkl operators for 
complex reflection groups introduced in \cite{DO}. Section 4 contains the main
result of our paper, Theorem 4.1, which states that the subalgebra of $\B(M_G)$ 
generated by the image of the map $\mu$ 
is isomorphic to the coinvariant algebra of a finite complex reflection group 
in question. Note that there is a duality 
between the corresponding NilCoxeter and coinvariant algebras 
in the case of finite real reflection groups, see e.g. 
\cite{Ba}. In Section 5 we study an analogue of such a duality for the group 
$G(e,1,n).$ Our results in Section 5 essentially depend on those obtained in 
\cite{RS}.  

\section{Coinvariant algebra of complex reflection group}
Let $G$ be a finite complex reflection group and $V$ the reflection 
representation of $G$. We fix a $G$-invariant hermitian inner product 
$\la \; , \; \ra$ on $V.$ Let $\A$ be the set of reflection hyperplanes 
$H\subset V$ of $G.$ 
The stabilizer $G_H \subset G$ of $H\in \A$ is isomorphic to a finite cyclic group 
$\Z/e_H \Z,$ $e_H \in \Z_{>0}.$ 
Each element $g\in G_H$ acts on $V$ as a complex reflection with 
respect to $H.$ We can assume that for $g\in G_H$ and $\xi \in V$ 
the action of $g$ is of form 
\[ g(\xi)= \xi -(1-\zeta)\frac{\la \xi , v_H \ra v_H}{\| v_H \|^2}, \]
where $v_H$ is a normal vector to $H$ and $\zeta$ is some 
$e_H$-th root of the unity. Denote by $\chi_H$ the character of the cyclic 
group $G_H$ defined as a restriction of $\det(g;V)$ to $G_H.$ 
For $H\in \A,$ there exists a unique 
element $g_H\in G_H$ such that $\chi_H (g)= \exp(2\pi \sqrt{-1}/e_H).$ 

Consider the symmetric algebra $S(V)={\rm Sym}_{\C}V$ 
of $V.$ 
The $G$-invariant 
subalgebra $S(V)^G$ is generated by algebraically independent homogeneous 
elements $f_1,\ldots , f_r,$ $r=\dim V,$ by Shephard and Todd's theorem 
\cite{ST}. The coinvariant algebra $P_G$ is the quotient algebra 
of $S(V)$ by the ideal $I_G$ generated by the fundamental $G$-invariants 
$f_1,\ldots , f_r.$ It has been shown by Chevalley \cite{Ch} that 
the algebra $P_G$ is isomorphic 
to the regular representation $\C \langle G \rangle$ as a left $G$-module. 

Let us fix a set of vectors $\{ v_H \}_{H\in \A }$ such that 
$v_H \in H^{\perp}\setminus \{ 0 \} .$ 
Then the action of $g\in G$ can be written as $g(v_H)=\lambda(g,H)v_{gH}$ for 
a constant $\lambda(g,H)\in \C^{\times}$ determined by $g$ and $H.$ The constants 
$\lambda(g,H)$ satisfy $\lambda({\rm id},H)=1$ and the cocycle condition 
$\lambda(gg',H)=\lambda(g,g'H)\lambda(g',H).$ 
The constants $\{ \lambda(g,H) \}_{g\in G,H\in \A}$ determine an element 
in $H^1(G,(\C^{\times})^{\A}).$ The family $\{ \alpha_H \}_{H\in \A}$ of 
defining linear forms of the reflection hyperplanes is also 
determined by $\alpha_H(x)=\la x, v_H \ra .$ Note that 
$g^* \alpha_H= \overline{\lambda(g^{-1},H)}\alpha_{g^{-1}H}$ for $g\in G.$ 

\begin{dfn} 
We define the divided difference operators $\Delta_{H,k}:S(V) \rightarrow 
S(V)$ as a $\C$-linear map defined by the formula 
\[ \Delta_{H,k}(f)=\frac{f-g_H^k(f)}{v_H} \] 
for $H\in \A$ and $1\leq k \leq e_H-1.$ 
\end{dfn}
The divided difference operators $\Delta_{H,k}$ satisfy the twisted 
Leibniz rule: 
\[ \Delta_{H,k}(f_1f_2)=\Delta_{H,k}(f_1)f_2+g_H^k(f_1)\Delta_{H,k}(f_2). \] 
It is easy to see the following. 
\begin{lem} 
A polynomial $f \in S(V)$ is a $G$-invariant if and only if 
$\Delta_{H,k}(f)=0$ for any $H\in \A$ and $1\leq k \leq e_H-1.$ 
\end{lem}
\section{Nichols-Woronowicz algebra over complex reflection group}
In this section we introduce the Nichols-Woronowicz algebra 
associated to a Yetter-Drinfeld module $M_G$ over the complex reflection 
group $G.$ In general, 
the Yetter-Drinfeld module over a finite group $\Gamma$ is defined as 
follows: 
\begin{dfn}
A vector space $M$ is called a Yetter-Drinfeld module over $\Gamma,$ if 
the following conditions are satisfied: \\ 
$(1)$ $M$ is a $\Gamma$-module, \\
$(2)$ $M$ is $\Gamma$-graded, i.e. $V=\bigoplus_{g\in \Gamma}M_g,$ where $M_g$ is a 
linear subspace of $M,$ \\ 
$(3)$ for $h\in \Gamma$ and $v\in M_g,$ $h(v)\in M_{hgh^{-1}}.$ 
\end{dfn} 
Note that the category $_{\Gamma}^{\Gamma}YD$ of the Yetter-Drinfeld modules over 
a fixed finite group $\Gamma$ is naturally braided by the braiding 
\[ \begin{array}{cccc} 
\psi_{M_1,M_2}: & M_1 \otimes M_2 & \rightarrow & M_2 \otimes M_1 \\ 
 & x \otimes y & \mapsto  & g(y) \otimes x , 
\end{array} \]
where $M_1,M_2\in \, ^{\Gamma}_{\Gamma}YD$ and $x \in (M_1)_g.$ 

Fix a set of normal vectors $v_H,$ $H\in \A.$ Then we can consider 
the corresponding constants $\lambda(g,H),$ $g\in G,$ $H\in \A.$ 
Let $M_G$ be a $\C$-vector space generated by the symbols $[H;k],$ $H\in \A$ 
and $1\leq k \leq e_H-1.$ 

\begin{dfn}
We define a structure 
of the Yetter-Drinfeld module over $G$ on the space $M_G$ as follows: \\ 
{\rm (i)} {\rm ($G$-action)} $g([H;k])= \overline{\lambda(g,H)}^{-1}[gH;k],$ \\ 
{\rm (ii)} {\rm ($G$-grading)} $\deg_G([H;k])=g_H^k.$ 
\end{dfn} 

\begin{lem} 
The $G$-action and $G$-grading defined above satisfy the condition 
for the Yetter-Drinfeld module, i.e., $\deg_G(h([H;k]))=hg_H^k h^{-1}$ for 
$h\in G.$ 
\end{lem}
{\it Proof.} 
Since 
\[ g(h([H;k]))=g(\overline{\lambda(h,H)}^{-1}[hH;k])=
\overline{\lambda(g,hH)}^{-1}\overline{\lambda(h,H)}^{-1}[ghH;k] \] 
\[ =\overline{\lambda(gh,H)}^{-1}[ghH;k]=(gh)([H;k]), \] 
the formula in (i) defines a $G$-action on $M_G.$ 
Let us check the condition $(2).$ From the definition of the $G$-action, 
we have 
\[ \deg_G(h([H;k]))=\deg_G(\overline{\lambda(h,H)}^{-1}[hH;k])=\deg_G([hH;k])=g_{hH}^k=hg_H^k h^{-1}. \] 

\begin{rem}
{\rm Our definition of the Yetter-Drinfeld module $M_G$ is analogous to 
the construction for the Coxeter group given in \cite[Section 5]{MS}. 
In the case of finite Coxeter groups, we can choose the constants $\lambda(g,H)$ 
to take the values $\pm 1$ by the normalization $\| v_H \|=1$ as in the construction 
of the Yetter-Drinfeld module $V_W$ used in \cite{Ba}. 
However, it is essential to specify the cocycle $\{ \lambda(g,H) \}$ in our case 
because of the appearance of the multiplication by some root of the unity.}
\end{rem}

For a braided vector space $M$ with a braiding $\Psi:M\otimes M \rightarrow 
M\otimes M,$ consider the linear endomorphism 
$\Psi_i$ on $M^{\otimes n}$ obtained by applying 
the braiding $\Psi: M \otimes M \rightarrow M \otimes M$ on the $i$-th and $(i+1)$-st 
components of $M^{\otimes n}.$ The endomorphisms $\Psi_i$ satisfy the braid relation 
$\Psi_{i+1} \Psi_i \Psi_{i+1}=\Psi_i \Psi_{i+1} \Psi_i.$ 
Denote by $s_i$ the simple transposition $(i,i+1)\in S_n.$ 
For any reduced expression
$w=s_{i_1}\cdots s_{i_l}\in S_n,$ the endomorphism 
$\Psi_w:=\Psi_{i_1}\cdots \Psi_{i_l}:M^{\otimes n} 
\rightarrow M^{\otimes n}$ is well-defined. 
The Woronowicz symmetrizer (\cite{Wo}) is given by 
\[ \sigma_n := \sum_{w\in S_n} \Psi_w . \]

\begin{dfn}
The Nichols-Woronowicz algebra associated to a braided vector space $M$ is 
\[ \B(M):= \bigoplus_{n\geq 0} M^{\otimes n}/{\rm Ker}(\sigma_n), \] 
where $\sigma_n : M^{\otimes n} \rightarrow M^{\otimes n}$ is the braided 
symmetrizer. 
\end{dfn}

The braided vector space $M$ naturally acts on $\B(M^*)$ from the right via 
the right braided derivations $\overleftarrow{D}_x,$ $x\in M.$ 
When $\Psi_{M,T(M^*)}^{-1}(\psi \otimes x)=\sum_ix_i\otimes \psi_i,$ 
denote by $\Psi_{M,T(M^*)}^{-1}(\psi \otimes \overleftarrow{D}_x)$ the 
operator $\phi \mapsto \sum_i(\phi \overleftarrow{D}_{x_i})\psi_i.$ 
The operators $\overleftarrow{D}_x$ are determined by the braided 
Leibniz rule 
\[ (\phi \psi)\overleftarrow{D}_x=\phi(\psi \overleftarrow{D}_x)+ 
\phi \Psi_{M,T(M^*)}^{-1}(\psi \otimes \overleftarrow{D}_x) , \] 
and the condition $\varphi \overleftarrow{D}_x=\varphi(x),$ $\varphi \in M^*,$ 
$x\in M,$ see \cite[2.5]{Ba}. 
In the subsequent construction, we identify the Yetter-Drinfeld module 
$M_G$ with its dual $M_G^*$ via 
the $G$-invariant symmetric inner product on $M_G$ 
given by $([H;k],[H';k'])=\delta_{H,H'}\delta_{k,k'}.$  
In our case, we have $\Psi_{M_G,T(M_G)}([H;k] \otimes \phi)= g_H^k(\phi)\otimes 
[H;k].$ Hence, the braided Leibniz rule can be written as 
\[ (\phi \psi)\overleftarrow{D}_{H,k}=\phi(\psi \overleftarrow{D}_{H,k})+ 
(\phi \overleftarrow{D}_{H,k})g_H^{-k}(\psi) . \] 
\begin{lem} {\rm (\cite[Criterion 3.2]{Ba})}
\[ \bigcap_{H\in \A, 1\leq k \leq e_H-1} {\rm Ker}(\overleftarrow{D}_{H,k})=\B^0(M_G) \] 
\end{lem}
The linear map 
\[ \begin{array}{cccc} 
\nu : & M_G & \rightarrow & {\rm End}_{\C}(\B(M_G)) \\  
 & [H;k] & \mapsto & \overleftarrow{D}_{H,e_H-k} 
\end{array} \]  
extends to the algebra homomorphism from the opposite algebra 
$\B(M_G)^{op}$ of $\B(M_G)$ to ${\rm End}_{\C}(\B(M_G)).$ 
The homomorphism $\nu:\B(M_G)^{op} \rightarrow {\rm End}_{\C}(\B(M_G))$ 
gives a nondegenerate pairing 
\[ \begin{array}{cccc} 
\lla \; , \; \rra : & \B(M_G) \times \B(M_G)^{op} & \rightarrow & \B^0(M_G)=\C \\ 
 & \lla \phi,\psi \rra & \mapsto & {\rm CT.}(\nu(\psi)(\phi)), 
\end{array} \] 
where ${\rm CT.}$ stands for the part of degree zero. 
\section{Relations in the Nichols-Woronowicz algebra} 
Denote by $d(N_1,N_2)>0$ the greatest common divisor of integers $N_1$ and 
$N_2.$ 
\begin{prop}
In the algebra $\B(M_G),$ 
\[ [H;k]^{e_H/d(e_H,k)}=0 . \] 
\end{prop}
{\it Proof.} Take a permutation $w\in S_n$ with $l(w)=l.$ Then 
\[ \Psi_w([H;k]^{\otimes n})=\zeta_H^{kl}[H;k]^{\otimes n} , \] 
where $\zeta_H=\exp(2\pi \sqrt{-1}/e_H).$ Hence 
\[ \sigma_n([H;k]^{\otimes n})= \left( \sum_{w\in S_n}\zeta_H^{kl(w)} \right)
\cdot [H;k]^{\otimes n} = 
\prod_{j=1}^{n-1}(1+\zeta_H^{k}+\zeta_H^{2k}+\cdots +\zeta_H^{jk})\cdot 
[H;k]^{\otimes n}. \] 
If $n=e_H/d(e_H,k),$ then $\sigma_n([H;k]^{\otimes n})=0.$ 
\bigskip \\ 
{\bf Relations for $G(e,1,n),$ $(e>1)$} \medskip \\ 
Take an $n$-dimensional hermitian vector space $V=\bigoplus_{i=1}^{n}\C \epsilon_i$ 
with an orthonormal basis $(\epsilon_i)_{i=1}^n.$ Let $(x_i)_{i=1}^n$ be the 
coordinate system with respect to the basis $(\epsilon_i)_{i=1}^n.$ 
All the reflection hyperplanes for $G(e,1,n) \subset GL(V)$ are given by 
\[ H_{ij}(a): x_i-\zeta^a x_j=0, \; \; \; H_i: x_i=0, \] 
where $1\leq i<j \leq n,$ $ a \in \Z / e \Z,$ $\zeta=\exp(2\pi \sqrt{-1}/e).$ 
Choose the normal vectors $v_{H_{ij}(a)}:=\epsilon_i - \zeta^{-a}\epsilon_j$ and 
$v_{H_i}:=\epsilon_i.$ 
The algebra $\B(M_G)$ is generated by the symbols $[H_{ij}(a)]:=[H_{ij}(a);1],$ 
$a\in \Z/e\Z,$ and 
$[H_i;s],$ $1\leq s \leq e-1.$ We put $[H_{ji}(a)]:=-\zeta^a [H_{ij}(-a)].$ 
\begin{prop} For $1\leq s,t \leq e-1$ and distinct $1\leq i,j,k \leq n,$ 
we have the following relations in $\B(M_G).$ \smallskip \\ 
$(1)$ $[H_{ij}(a)][H_{jk}(b)]- \zeta^{a}[H_{ik}(a+b)][H_{ij}(a)] - 
[H_{jk}(b)][H_{ik}(a+b)] =0,$ \smallskip \\
$(2)$ $\sum_{p=1}^{e_H/d(e_H,2(a-b))} \zeta^{2p(a-b)}[H_{ij}(a+2p(a-b))][H_{ij}(b+2p(a-b))]=$ \smallskip \\  
$\sum_{q=1}^{e_H/d(e_H,2(a-b))}\zeta^{(2q-1)(a-b)}[H_{ij}(b+2q(a-b))][H_{ij}(a+2(q-1)(a-b))],$ \smallskip \\
$(3)$ $[H_{ij}(a)][H_i;s]-\zeta^{-a}[H_j;s][H_{ij}(a)]+\zeta^{-a}[H_{ij}(a-s)][H_j;s]-\zeta^{-s}[H_i;s][H_{ij}(a-s)]=0,$ \smallskip \\   
$(4)$ $[H_{ij}(a)][H_{kl}(b)]=[H_{kl}(b)][H_{ij}(a)],$ if $\{ i,j \} \cap \{ k,l \} =\emptyset,$ \smallskip \\ 
$(5)$ $[H_i;s][H_j;t]=[H_j;t][H_i;s],$ \smallskip \\
$(6)$ $[H_{ij}(a)][H_k;s]=[H_k;s][H_{ij}(a)],$ if $k\not=i,j.$ 
\end{prop}
These relations follow from straightforward computation of the images 
of the braided symmetrizers. 
\begin{rem} 
{\rm In the case of $S_n=G(1,1,n),$ the relations $(1),$ $(4)$ and 
$[H_{ij}(0)]^2=0$ cover all the independent quadratic 
relations in $\B(M_{S_n}),$ see \cite{FK} and \cite[7.1]{Ba}. 
In the case of the Weyl groups $W_{B_n}=W_{C_n}=G(2,1,n),$ the relations 
in Proposition 3.2 coincide with the quadratic relations given in \cite{KM1}.}
\end{rem}  
\section{Model of coinvariant algebra}
\begin{dfn}
Fix a set of $G$-invariant constants $\ka=(\ka_{H,i})_{H\in \A,1\leq i \leq e_H-1}.$ 
We define the $\C$-linear map $\mu=\mu_{\ka}:V \rightarrow M_G$ by 
\[ \mu(\xi)= -\sum_{H\in \A}\sum_{i,k=1}^{e_H-1} \alpha_H(\xi) \ka_{H,i} \zeta_H^{-ik} [H;k], \] 
where $\zeta_H=\exp(2\pi \sqrt{-1}/e_H).$ 
\end{dfn}
\begin{prop}
The map $\mu$ is a $G$-homomorphism, i.e., $\mu(g(\xi))=g(\mu(\xi)).$ 
\end{prop}
{\it Proof.} Note that $e_{gH}=e_H$ for any $g\in G$ and $H\in \A.$ We have 
\begin{eqnarray*}
\mu(g(\xi))& = & -\sum_{H\in \A}\sum_{i,k=1}^{e_H-1} \alpha_H(g(\xi)) \ka_{H,i}\zeta_H^{-ik} [H;k] \\
 & = & -\sum_{H\in \A}\sum_{i,k=1}^{e_H-1} g^*\alpha_H(\xi) \ka_{H,i} \zeta_H^{-ik} [H;k] \\ 
 & = & -\sum_{H\in \A}\sum_{i,k=1}^{e_H-1} \overline{\lambda(g^{-1},H)}\alpha_{g^{-1}H}(\xi)  \ka_{H,i} \zeta_H^{-ik} [H;k] \\ 
 & = & -\sum_{H\in \A}\sum_{i,k=1}^{e_H-1} \alpha_H(\xi)  \ka_{gH,i} \zeta_{gH}^{-ik} \overline{\lambda(g^{-1},gH)}[gH;k] \\ 
 & = & -\sum_{H\in \A}\sum_{i,k=1}^{e_H-1} \alpha_H(\xi)  \ka_{H,i} \zeta_{H}^{-ik} \overline{\lambda(g,H)}^{-1}[gH;k] \\ 
 & = & g\left( -\sum_{H\in \A}\sum_{i,k=1}^{e_H-1} \alpha_H(\xi) \ka_{H,i} \zeta_{H}^{-ik} [H;k]\right) =g(\mu(\xi)) .  
\end{eqnarray*}
Here we have used that $\lambda(g^{-1},gH)=\lambda(g,H)^{-1}.$ 
\begin{prop} {\rm 
\[ [\mu(\xi),\mu(\eta)]=0 \; \; \; \textrm{in $\B(M_G)$} \] }
\end{prop}
{\it Proof.} Let us show $({\rm id}+\Psi)(\mu(\xi) \otimes \mu(\eta))=
({\rm id}+\Psi)(\mu(\eta) \otimes \mu(\xi)).$ 
The left-hand side equals 
\[ \sum_{H,H'\in \A}\sum_{i,k=1}^{e_H-1}\sum_{j,l=1}^{e_{H'}-1} \alpha_{H}(\xi)\alpha_{H'}(\eta) 
\ka_{H,i}\ka_{H',j}\zeta_H^{-ik}\zeta_{H'}^{-jl} [H;k]\otimes[H';l] \] 
\[ + \sum_{H,H'\in \A}\sum_{i,k=1}^{e_H-1}\sum_{j,l=1}^{e_{H'}-1} \alpha_{H}(\xi)\alpha_{H'}(\eta) 
\ka_{H,i}\ka_{H',j} \zeta_H^{-ik}\zeta_{H'}^{-jl} g_H^k([H';l])\otimes[H;k] . \] 
Here the second term is 
\begin{eqnarray*}  
 & & \sum_{H,H'\in \A}\sum_{i,k=1}^{e_H-1}\sum_{j,l=1}^{e_{H'}-1} 
\alpha_{H}(\xi)\alpha_{H'}(\eta)
\ka_{H,i}\ka_{H',j}\zeta_H^{-ik}\zeta_{H'}^{-jl} \overline{\lambda(g_H^k,H')}^{-1}[g_H^k(H');l]\otimes[H;k] \\  
 & = & \sum_{H,H'\in \A}\sum_{i,k=1}^{e_H-1}\sum_{j,l=1}^{e_{H'}-1} \alpha_H(\xi) 
\overline{\lambda(g_H^{-k},g_H^k(H'))}\alpha_{H'}(\eta ) \ka_{H,i}\ka_{H',j}
\zeta_H^{-ik}\zeta_{H'}^{-jl} [g_H^k(H');l]\otimes[H;k]  \\ 
 & = & \sum_{H,H'\in \A}\sum_{i,k=1}^{e_H-1}\sum_{j,l=1}^{e_{H'}-1} \alpha_H(\xi)
(g_H^{k})^*(\alpha_{g_H^k(H')})(\eta) \ka_{H,i}\ka_{H',j}
\zeta_H^{-ik}\zeta_{H'}^{-jl} [g_H^k(H');l]\otimes[H;k]  \\ 
 & = & \sum_{H,H'\in \A}\sum_{i,k=1}^{e_H-1}\sum_{j,l=1}^{e_{H'}-1} \alpha_H(\xi)  
\alpha_{H'} (g_H^{k}(\eta)) \ka_{H,i}\ka_{H',j}
\zeta_H^{-ik}\zeta_{H'}^{-jl} [H';l]\otimes[H;k] . 
\end{eqnarray*} 
Since 
\[ g_H^k(\eta)=\eta-(1-\zeta_H^{k})\frac{\alpha_H (\eta) v_H}{ \| v_H \|^2}, \] 
we obtain the following expression of $({\rm id}+\Psi)(\mu(\xi) \otimes \mu(\eta))$ 
which is symmetric in $\xi$ and $\eta:$ 
\[ \sum_{H,H'\in \A}\sum_{i,k=1}^{e_H-1}\sum_{j,l=1}^{e_{H'}-1} \alpha_{H}(\xi)\alpha_{H'}(\eta) 
\ka_{H,i}\ka_{H',j} \zeta_H^{-ik}\zeta_{H'}^{-jl} [H;k]\otimes[H';l] \] 
\[ + \sum_{H,H'\in \A}\sum_{i,k=1}^{e_H-1}\sum_{j,l=1}^{e_{H'}-1} \alpha_{H'}(\xi) \alpha_{H}(\eta) \ka_{H,i}\ka_{H',j}
\zeta_H^{-ik}\zeta_{H'}^{-jl} [H;k]\otimes[H';l] \]
\[ + \sum_{H,H'\in \A}\sum_{i,k=1}^{e_H-1}\sum_{j,l=1}^{e_{H'}-1} \alpha_{H'}(\xi) \alpha_{H'}(\eta)
\frac{\alpha_{H}(v_{H'})}{ \|v_{H'} \|^2} \ka_{H,i}\ka_{H',j}
\zeta_H^{-ik}\zeta_{H'}^{-jl}(1-\zeta_{H'}^{l}) [H;k]\otimes[H';l] . \]
This completes the proof. \bigskip 

The proposition above shows that the map $\mu$ extends to an algebra homomorphism 
\[ \tilde{\mu} : S(V) \rightarrow \B(M_G) . \] 
\begin{rem}{\rm 
Dunkl and Opdam \cite{DO} have introduced the Dunkl operator for the complex 
reflection groups. Their operator is defined by the following formulas: 
\begin{eqnarray*}  
T_{\xi}(\ka) & = & \partial_{\xi}+\sum_{H\in \A}\sum_{i=1}^{e_H-1}\sum_{g\in G_H}
\alpha_H(\xi)\ka_{H,i}\alpha_H^{-1} \chi^i_H(g)g \\ 
 & = & \partial_{\xi}-\sum_{H\in \A}\sum_{i=1}^{e_H-1}\sum_{k=1}^{e_H-1}
\alpha_H(\xi)\ka_{H,i} \zeta_H^{-ik} \alpha_H^{-1}(1-g_H^{-k}), 
\end{eqnarray*}
where $\chi^i_H$ is the restriction of $\det^i$ to $G_H.$ 
Hence, 
our homomorphism $\mu$ can be regarded as a truncated version of the operator 
$T_{\xi}(\ka)$ that means the operator without the differential part $\partial_{\xi}$ after replacing 
the brackets $[H;k]$ by the divided difference operators $\alpha_H^{-1}(1-g_H^{-k})$ 
acting on $S(V^*).$}
\end{rem}

Below we quote a lemma from \cite{RS}, which is an analogue 
of Lemma 1.9 in \cite[Chapter IV]{Hi2}. The proof of this lemma given 
in \cite{RS} is applicable to the general finite complex reflection 
groups. Let us consider 
the polynomial $Q:=\prod_{H\in \A} v_H^{e_H-1} \in S(V).$ 
\begin{lem} {\rm (\cite[Lemma 2.16]{RS})} 
If a graded ideal $I\subset S(V)$ contains $I_G,$ but does not 
contain $Q,$ then $I=I_G.$ 
\end{lem}

We define the divided difference operator $\overleftarrow{\Delta}_{H,k}$ 
acting from the right by 
\[ f\overleftarrow{\Delta}_{H,k}:=\frac{f-g_H^{-k}(f)}{v_H}. \]
\begin{prop} For $f\in S(V),$ 
\[ C_{H,k}\tilde{\mu}(f\overleftarrow{\Delta}_{H,k})=
\tilde{\mu}(f)\overleftarrow{D}_{H,k} , \] 
where 
\[ C_{H,k}=\sum_{i=1}^{e_H-1}\ka_{H,i}\frac{\zeta_H^{-ik}\| v_H \|^2}{\zeta_H^{-k}-1} . \] 
\end{prop}
{\it Proof.} This follows from 
\[ [H';k']~\overleftarrow{D}_{H,k}=\delta_{H,H'}\delta_{k,k'}, \] 
\[ \mu(\xi)\overleftarrow{D}_{H,k}= 
(-\sum_{i=1}^{e_H-1}\ka_{H,i}\zeta_H^{-ik}) \alpha_H(\xi) = 
\left( \sum_{i=1}^{e_H-1}\ka_{H,i}\frac{\zeta_H^{-ik}\| v_H \|^2}{\zeta_H^{-k}-1}\right) 
\xi\overleftarrow{\Delta}_{H,k}, 
\; \; \; \textrm{for $\xi \in V,$} \] 
and the braided Leibniz rule 
\[ (\phi \psi)\overleftarrow{D}_{H,k}=\phi (\psi \overleftarrow{D}_{H,k})+
(\psi \overleftarrow{D}_{H,k})g_H^{-k}(\phi). \] 
\medskip 

The constants $(\ka_{H,i})$ are said to be {\it generic} if 
$C_{H,k}\not= 0$ for any $H\in \A$ and $1\leq k\leq e_H-1.$ 
We also need the following lemma to prove the main theorem. 
\begin{lem} Let $(\ka_{H,i})$ be generic. 
There exist sequences $H_1,\ldots,H_p\in \A$ and $k_1,\ldots, k_p$ such that 
$Q\overleftarrow{\Delta}_{H_1,k_1}\cdots \overleftarrow{\Delta}_{H_p,k_p}$ is a 
nonzero constant. 
\end{lem}
{\it Proof.} Take a nonzero homogeneous element $f\in P_G.$ 
If $f\overleftarrow{\Delta}_{H,k}=0$ for any $H$ and $k,$ then 
$f$ is a $G$-invariant element of $P_G$ by Lemma 1.1. 
Since $P_G$ is isomorphic to the 
regular representation of $G$ and contains only one copy of the trivial 
representation at degree zero, $f$ must be in $P_G^0=\C.$ Hence, if 
the degree of $f$ is positive, then there exist $H\in \A$ and $k$ such 
that $f\overleftarrow{\Delta}_{H,k}\not=0.$ 
The polynomial $Q$ affording the character $\det^{-1}$ 
is a generator of the part of the highest degree in $P_G,$ 
so we can find the desired sequences 
$H_1,\ldots,H_p\in \A$ and $k_1,\ldots, k_p$ by induction on the degree. 

\begin{thm} For any generic choice of the constants $\ka=(\ka_{H,i}),$ we have 
the isomorphism 
\[ {\rm Im}(\tilde{\mu}) \cong P_G. \] 
\end{thm}
{\it Proof.} If $f$ is a $G$-invariant polynomial 
of positive degree, we get $\tilde{\mu}(f)\overleftarrow{D}_{H,k}=0$ 
for any $H$ and $k$ from Lemma 1.2 and Proposition 4.3. Hence 
$\tilde{\mu}(f)\in \B^0(M_G)=\C$ from Lemma 2.2. 
Since $\tilde{\mu}$ preserves the degree, we have $\tilde{\mu}(f)=0$ 
and ${\rm Ker}(\tilde{\mu}) \supset I_G.$ On the other hand, it follows 
from Lemma 4.2 that 
one can find sequences 
$H_1,\ldots,H_p\in \A$ and $k_1,\ldots, k_p$ such that 
$Q\overleftarrow{\Delta}_{H_1,k_1}\cdots \overleftarrow{\Delta}_{H_p,k_p}$ is a 
nonzero constant. This means, see Proposition 4.3, that 
\[ \tilde{\mu}(Q)\overleftarrow{D}_{H_1,k_1}\cdots \overleftarrow{D}_{H_p,k_p}
= C_{H_1,k_1}\cdots C_{H_p,k_p} \tilde{\mu}(Q\overleftarrow{\Delta}_{H_1,k_1}\cdots 
\overleftarrow{\Delta}_{H_p,k_p})  \in \B^0(M_G) \] 
is a nonzero constant. Hence ${\rm Ker}(\tilde{\mu})$ does not contain $Q.$ 
Now we get ${\rm Ker}(\tilde{\mu})=I_G$ from Lemma 4.1.  

\section{Complex reflection group of type $G(e,1,n)$}
We use the notation in Section 3. 
In the following we identify $V^*$ with $V$ via the $G$-invariant hermitian form $\la \; , \; \ra,$ 
so that the left $G$-action on $V^*$ defined by $g(\alpha)=(g^{-1})^*\alpha,$ 
$\alpha \in V^*,$ coincides with the left $G$ action on $V.$  
Put $\ka_{ij,a}=\ka_{H_{ij}(a),1}$ and 
$\ka_{i,s}=\ka_{H_i,s}.$ 
For a generic choice of the constants $\ka,$ the subalgebra ${\rm Im}(\tilde{\mu}) \subset \B(M_G),$ 
which is generated by the elements 
\[ \mu(\epsilon_i)= \sum_{j\not=i}\sum_{a\in \Z/m\Z}\ka_{ij,a}[H_{ij}(a)] - 
\sum_{s,t=1}^{e-1}\zeta^{-st}\ka_{i,s}[H_i;s], \; \; 1\leq i \leq n, \] 
is isomorphic to the coinvariant algebra 
\[ P_G= \C[x_1,\ldots,x_n]/(E_1(x_1^e,\ldots,x_n^e),\ldots,E_n(x_1^e,\ldots,x_n^e)), \] 
where $E_i$ is the $i$-th elementary symmetric polynomial. 
Rampetas and Shoji \cite{RS} introduced a family of operators $\Delta_w,$ $w\in G,$ 
acting on the polynomial ring $P=\C[x_1,\ldots,x_n]$ based on particular 
choice of the reduced expression of $w\in G.$ The pseudo-reflections 
$s_{ij}(a):=g_{H_{ij}(a)}$ and $t_i:=g_{H_i}$ generate the group $G(e,1,n).$ 
In particular, the pseudo-reflections $s_i=s_{i-1 \, i}(0)$ and $t_1$ play 
the role of simple reflections. 
In the following, we use the divided difference operators 
$\Delta_{s_i}:=\Delta_{H_{i-1\, i}(0),1}$ and $\Delta_{t_i}:=\Delta_{H_i,1}.$ 
Note that the braid relations 
among the divided difference operators do not hold in general. 
Put $\tilde{\Delta}_{s_i}:=\Delta_{t_i}^{e-2}\Delta_{s_i}.$ 
Let us consider the operator 
\[ \Delta_n(k,a)= \left\{ \begin{array}{cc} 
\Delta_{s_{k+1}} \cdots \Delta_{s_{n-1}}\Delta_{s_n}, & \textrm{if $a=0,$} \\ 
\tilde{\Delta}_{s_k}\cdots \tilde{\Delta}_{s_2} \Delta_{t_1}^a 
\Delta_{s_2} \cdots \Delta_{s_n}, & \textrm{if $a\not=0.$} 
\end{array} \right. \]
In \cite{RS}, it is shown that any element $w\in G$ has a unique 
decomposition of form 
\[ w=\omega_n(k_n,a_n)\omega_{n-1}(k_{n-1},a_{n-1})\cdots \omega_1(k_1,a_1), \] 
where 
\[ \omega_n(k,a)= \left\{ \begin{array}{cc} 
s_{k+1} \cdots s_{n-1}s_n, & \textrm{if $a=0,$} \\ 
s_k \cdots s_2 t_1^a s_2 \cdots s_n, & \textrm{if $a\not=0.$} 
\end{array} \right. \]
The operators $\Delta_w$ are defined by the following formula: 
\[ \Delta_w:= \Delta_n(k_n,a_n)\Delta_{n-1}(k_{n-1},a_{n-1})\cdots \Delta_1(k_1,a_1). \] 
Define the evaluation map $\varepsilon:P \rightarrow \C$ by $\varepsilon(f)=f(0).$ 
Denote by $\bar{\cal{D}}_G$ the subspace of $P^*$ spanned by 
the operators $\varepsilon \Delta_w,$ $w\in G.$ 
\begin{prop} {\rm (\cite[Theorem 2.18]{RS})} 
The coinvariant algebra $P_G$ is naturally isomorphic to the dual 
space of $\bar{\cal{D}}_G.$ 
\end{prop}

Put $[s_i]:=[H_{i-1\, i}(0)],$ $[t_i]:=[H_i;e-1]$ and $\widetilde{[s_i]}:=[t_i]^{e-2}[s_i].$ 
We also define the elements $[w] \in \B(M_G)^{op},$ $w\in G,$ 
in a similar way: 
\[ [w]:=[\omega_n(k_n,a_n)][\omega_{n-1}(k_{n-1},a_{n-1})]\cdots [\omega_1(k_1,a_1)], \] 
where 
\[ [ \omega_n(k,a) ] := \left\{ \begin{array}{cc} 
[ s_{k+1} ] \cdots [s_{n-1}][ s_n ], & \textrm{if $a=0,$} \\ 
\widetilde{[s_k]} \cdots \widetilde{[s_2]} [t_1]^a [s_2] \cdots [s_n], & \textrm{if $a\not=0.$} 
\end{array} \right. \] 
Consider ${\bf D}_G$ the subspace of $\B(M_G)^{op}$ spanned by the elements 
$[w],$ $w\in G.$ Now it is easy to get the following analogue to \cite[Theorem 6.1]{Ba}. 
\begin{thm} Assume that $\ka_{ij,a}=1$ and $\ka_{i,s}=1-\zeta^{-s}.$ \smallskip \\ 
$(1)$ The linear map 
\[ ({\rm CT.})_*\circ \tilde{\mu}^* \circ \nu: \B(M_G)^{op} \rightarrow {\rm End}_{\C}(\B(M_G)) \rightarrow 
{\rm Hom}_{\C}(P,\B(M_G)) \rightarrow P^* \] 
induces an isomorphism between 
${\bf D}_G$ and $\bar{\cal{D}}_G.$ \\ 
$(2)$ The subalgebra ${\rm Im}(\tilde{\mu})$ is isomorphic to the dual 
space of ${\bf D}_G$ via the pairing $\lla \; , \; \rra.$ Furthermore, 
the pairing $\lla \; , \; \rra$ restricted to ${\rm Im}(\tilde{\mu}) \times 
{\bf D}_G$ coincides with the pairing between $P_G$ and $\bar{\cal{D}}_G,$ i.e., 
\[ \lla \tilde{\mu}(f) , [w] \rra = \varepsilon \Delta_w(f). \] 
\end{thm} 
{\it Proof.} From Theorem 4.1 and Proposition 5.1, we have the factorization  
\[ ({\rm CT.})_* \circ \tilde{\mu}^*: {\rm End}_{\C}(\B(M_G)) \rightarrow P_G^*=\bar{\cal{D}}_G \subset P^*. \] 
For the choice of the constants $\ka$ as assumed, 
we have $C_{H,k}=1$ for all $H\in \A$ and $1 \leq k \leq e_H-1.$ 
Proposition 4.3 shows that 
\[ \lla \tilde{\mu}(f),[w] \rra = {\rm CT.}(\nu([w])(\tilde{\mu}(f)))= {\rm CT.}(\tilde{\mu}(\Delta_w(f)))=
\varepsilon \Delta_w(f), \; \; f\in P, \] 
so we obtain $({\rm CT.})_*\tilde{\mu}^*(\nu({\bf D}_G))=\bar{\cal{D}}_G.$ 
Since ${\bf D}_G$ is spanned by $|G|$ elements, ${\bf D}_G$ is isomorphic to 
$\bar{\cal{D}}_G.$ 
 
Anatol N. Kirillov \\
Research Institute for Mathematical Sciences \\ 
Kyoto University \\
Sakyo-ku, Kyoto 606-8502, Japan \\
e-mail: {\tt kirillov@kurims.kyoto-u.ac.jp} \\
URL: {\tt http://www.kurims.kyoto-u.ac.jp/\textasciitilde kirillov} 
\bigskip \\ 
Toshiaki Maeno \\
Department of Electrical Engineering, \\
Kyoto University, \\ 
Sakyo-ku, Kyoto 606-8501, Japan \\ 
e-mail: {\tt maeno@kuee.kyoto-u.ac.jp}

\end{document}